\documentclass{amsart}

\usepackage{amsfonts}
\usepackage{amssymb}
\usepackage{verbatim}

\author{Caleb D. Holloway and Tavan T. Trent}
\title{Wolff's Theorem on Ideals for Matrices}

\newtheorem{mydef}{Definition}[section]

\newtheorem{lem}{Lemma}[section]

\newtheorem*{wolffmatrix}{Wolff's Theorem for Matrices}

\newcommand{\ran}{\operatorname{ran}}
\newcommand{\rank}{\operatorname{rank}}
\newcommand{\all}{\; \forall \;}

\newcommand{\C}{\mathbb{C}}
\newcommand{\F}{\mathcal{F}}
\newcommand{\G}{\mathcal{G}}
\newcommand{\A}{\mathcal{A}}
\newcommand{\I}{\mathcal{I}}
\newcommand{\D}{\mathbb{D}}
\newcommand{\Dir}{\mathcal{D}}
\newcommand{\sgn}{\text{sgn}}

\allowdisplaybreaks[1]

\begin{document}

\maketitle

\begin{abstract}
We extend Wolff's theorem concerning ideals on $ H^{\infty}(\D) $ to the matrix case, giving conditions under which an $ H^{\infty} $-solution $ G $
to the equation $ FG = H $ exists for all $ z \in \D $, where $ F $ is an $ m \times \infty $ matrix of functions in $ H^{\infty}(\D) $, and $ H $ is an
$ m \times 1 $ vector of such functions. We then examine several useful results.
\end{abstract}

\section{Introduction}

In 1962, Lennart Carleson \cite{carleson} solved the Corona Problem, proving that the ideal $ \I $ generated by a finite set of functions $ \lbrace f_i \rbrace_{i=1}^{n} \subset H^{\infty}(\D) $ is the entire space $ H^{\infty}(\D) $ provided there exists $ \delta > 0 $ such that
\begin{equation}
\sum_{i=1}^{n} | f_i(z) | \geq \delta \all z \in \D \text{.}
\end{equation}
This result can be extended to hold for infinitely many functions $ \lbrace f_i \rbrace_{i=1}^{\infty} $ (see \cite{rosenblum}, \cite{operatortheory}). Two different extensions of the corona theorem to matrices were given by Fuhrmann \cite{fuhrmann} and Andersson \cite{andersson}. Fuhrmann's result 
was extended to one-sided infinite matrices by Vasyunin (see Nikolski \cite{nikolski}). 
However, Treil \cite{treil} showed that a complete extension of the corona 
theorem is not possible in the two-sided infinite matrix case. Trent and 
Zhang proved that the result of Furmann and Vasyunin can be 
extended to any algebra that satisfies a corona theorem \cite{mct}, and later did the same for
Andersson's result, also allowing for one-sided infinite matrices \cite{mctII}.

A more general question than Carleson's is, under what conditions is a given function $ h \in H^{\infty}(\D) $ to be found in $ \I $? One might suppose based on Carleson's result that a sufficient condition would be that 
\begin{equation}\label{eq:wolff}
\sum_{i=1}^{n} |f_i(z)| \geq |h(z)| \all z \in \D \text{,}
\end{equation}
but that is not the case, as Rao proved (see Garnett \cite{garnett}). Thomas Wolff, however, proved that, given \eqref{eq:wolff}, $ h^3 \in \I $ \cite{wolff1}. More recently, Treil showed that the result fails when the exponent ``3'' is replaced with ``2'' (although it holds for any exponent greater than 2) \cite{treil2}.

If we adjust the hypothesis to 
\begin{equation} \label{eq:scalar}
[\; \sum_{i=1}^{n} |f_i(z)|^2 \;]^{\frac{3}{2}} \geq |h(z)| \all z \in \D \text{,}
\end{equation}
we obtain $ h \in \I $. A matter of current interest is how this estimate can be improved, which we will discuss at the end of this paper. For now,  \eqref{eq:scalar} will be sufficient for our needs.

Since the corona theorem has been applied to matrices, a natural question is whether the same can be done for Wolff's theorem. The answer is yes, as we will now set forth to prove. First, however, a definition is in order.

\begin{mydef}
Let $ B \in M_m(\C) $. For $ 1 \leq k \leq m $, define
\begin{equation*}
{\det}_k(B) = \sum \limits_{\pi \in \Pi_k(m)} \det(E_{\pi} B E_{\pi})
\end{equation*}
where $ \Pi_k(m) $ denotes the increasing k-tuples of integers in $ \lbrace 
1, 2, \dots , m \rbrace $.
\end{mydef}

Here $ E_{\pi} $ is the $ m \times m $ matrix whose $ i $th column is the $ 
i $th column of the $ m \times m $ identity matrix if $ i \in \pi $, and is 
zero otherwise. When taking the determinant of $ E_{\pi} B E_{\pi} $ in the 
above definition, we delete those columns and rows consisting of all zeros.

\begin{wolffmatrix} \label{wolffmatrix}
Let $ F(z) $ be an $ m \times \infty $ matrix of functions in $ 
H^{\infty}(\D) $ with $ \max\{\rank  F(z) \mid z \in \D \} = k \leq m $.
Let $ H(z) $ be an $ m \times 1 $ vector of functions in $ H^{\infty}(\D) 
$. Suppose
\begin{itemize}
\item[(i)] $ [ \det_k (F(z)F(z)^{\ast})]^{\frac{3}{2}} \geq |h_i(z)| \all  
z \in \D $, $ i = 1, \dots , m $
\item[(ii)] $ \| M_F \| = 1 $
\item[(iii)] there exists a function $ \underline{u} : \D \rightarrow l^2 $ 
such that $ F \underline{u} = H $ everywhere on $ \D $.
\end{itemize}
Then there exists an $ \infty \times 1 $ vector $ G(z) $ of functions in $ 
H^{\infty}(\D) $ such that
\begin{itemize}
\item[(a)] $ F(z)G(z) = H(z) \all z \in \D, $ and
\item[(b)] $ \| M_G \| < \infty $.
\end{itemize}
\end{wolffmatrix}

We base our arguments on those found in Trent and Zhang \cite{mctII}. The main difference here is that we do not assume a uniform lower boundedness on $ F $, and instead assume that $ F $ is bounded by the entries of $ H $.

\section{Preliminaries}

Before giving the proof of Wolff's Theorem for Matrices, we define and list some properties of ``Q-operators.'' Proofs of these properties can be found in \cite{mct}.

We let $ H \wedge K $ denote the exterior product between two Hilbert spaces $ H $ and $ K $, and $ l^2_{(n)} = \wedge_{i=1}^{n}l^2 $. In keeping with this notation, $ l^2_{(0)} = \C $.

Let $ \lbrace e_i \rbrace_{i = 1}^{\infty} $ denote the standard basis in $ l^2 $. If $ I_n $ denotes increasing $n$-tuples of positive integers and if 
$ (i_1, i_2, \dots , i_n) \in I_n $, we let $ \pi_n = \lbrace i_1, i_2, 
\dots , i_n \rbrace $ and, abusing notation, we write $ \pi \in I_n $. If 
we define $ e_{\pi_n} = e_{i_1} \wedge e_{i_2} \wedge \dots \wedge e_{i_n} 
$, then $ \lbrace e_{\pi_n} \rbrace_{\pi_n \in I_n} $ is defined to be the 
standard basis 
for $ l^2_{(n)} $.

Let $ H(E) $ be a reproducing kernel Hilbert space on a set $ E $, and let 
$ \A = M(H(E)) $, the multiplier algebra on $ H(E) $. Let $ \underline{f}_j(z) = 
(v_1(z), v_2(z), \dots ) $, where $ \lbrace v_n \rbrace_{n = 1}^{\infty} 
\subset \A $, such that $ \underline{f}_j(z)\underline{f}_j(z)^{\ast} \leq 1 \all z \in E $. Fix $ 
z \in E $, and for $ n = 0, 1, \dots $ define
\begin{equation*}
Q_j^{(n)\ast}(z) : l^2_{(n)} \rightarrow l^2_{(n + 1)}
\end{equation*}
by
\begin{equation*}
Q_j^{(n)\ast}(z)(w_n) = \overline{\underline{f}_j(z)} \wedge \underline{w}_n \text{,}
\end{equation*}
where $ \underline{w}_n \in l^2_{(n)} $.  Note that $ Q_j^{(0)\ast}(z) = \overline{\underline{f}_j(z)} $.

We observe that $ \ran Q_j^{(n)\ast}(z) \subset \ker 
Q_j^{(n+1)\ast}(z) $. Furthermore, equality can be shown 
if we stipulate that $ \underline{f}_j(z)\underline{f}_j(z)^{\ast} \geq \delta > 0 $.\ (This follows
from \eqref{Qid} below.)

By the anti-commutivity of the exterior product, we see that
\begin{equation}
\label{anticommute}
Q_j^{(n)}Q_k^{(n+1)} = -Q_k^{(n)}Q_j^{(n+1)} \text{.}
\end{equation}
Also, for $ e_{\pi_n} \in l^2_{(n)} $, we have
\begin{equation*}
Q_j^{(n)\ast}(z)(e_{\pi_n}) = \overline{\underline{f}_j(z)} \wedge e_{\pi_n} = (\sum 
\limits_{p = 1}^{\infty} \overline{v_p(z)}e_p) \wedge e_{\pi_n} \text{,}
\end{equation*}
so with respect to the standard basis, the entries in $ Q_j^{(n)\ast}(z) $ 
are 0 or $ \pm \overline{v_n(z)} $ for some $ n $. Thus $ Q_j^{(n)}(\cdot) 
$ has entries belonging to $ \A $ with respect to the standard basis.

Assume that there exists a fixed $ z \in E $ such that $ f_j(z) \neq 0 $, 
and let $ Q_n = Q_j^{(n)}(z) $. Then for $ n = 0, 1, \dots $,
\begin{equation} \label{Qid}
Q_n^{\ast}Q_n + Q_{n+1}Q_{n+1}^{\ast} = \|a\|^2 I_{l^2_{n+1}} \text{.}
\end{equation}

Finally, if $ a_i = \underline{f}_i(z) $, then for $ k = 2, \dots , m $, let $ a_k' = 
P_{\text{sp} \lbrace a_1, \dots , a_{k-1} \rbrace}^{\perp}(a_k) $. Then
\begin{align}
\label{det}
a_1 Q_{a_2}^{(1)} \dots Q_{a_m}^{(m-1)} Q_{a_m}^{(m-1)\ast} \dots 
Q_{a_2}^{(1)\ast} a_1^{\ast} &= \|a_1\|^2 \prod_{j=2}^{m} \|a_j'\|^2 
\notag \\ &= \det(F(z)F(z)^{\ast}) \text{.}
\end{align}
The last equality is obtained by a straightforward computation, using the 
fact that $ \frac{a_i^{\ast}a_i}{\|a_i\|^2} $ is the rank one projection of 
$ l^2 $ onto $ a_i $.

We obtain a very nice norm estimate in the case where $ \A = 
H^{\infty}(\D) $. Let $ A = (a_1, a_2, \dots), a_i \in H^{\infty}(\D) $, with $ \|M_A\| < \infty $. From \eqref{Qid} we have the following estimate on the 
operator norms (for fixed $ z \in \D $):
\begin{equation*}
 \| Q^{(j)} \| \leq \| A \| \text{.}
\end{equation*}
Now
\begin{equation*}
\| M_{Q^{(j)}} \| = \sup_{z \in \D} \| Q^{(j)}(z) \| \leq \sup_{z \in \D} 
\| A(z) \| = \| M_A \| \text{.}
\end{equation*}
(Here $ j \geq 1 $. If $ j = 0 $, then $ Q^{(j)} = A $, so the result is 
trivial.) Thus
\begin{equation} \label{Hinfinitynormestimate}
\| M_{Q^{(j)}} \| \leq \| M_A \| \text{.}
\end{equation}
We will need Lemma 1 from \cite{mctII} if we wish to extend our main theorem to 
spaces besides $ H^{\infty}(\D) $, however.

The following lemmas will be used in the proof of our theorem. Their proofs have been omitted but may be found in \cite{mctII}.

\begin{mydef}
Let $ \lbrace X_j \rbrace_{j=1}^{n+1} $ be Banach spaces and let $ \lbrace 
T_{jk} \rbrace_{j,k = 1}^n $ denote operators such that $ T_{jk} \in 
B(X_{j+1}, X_j) $. Let
\begin{equation*}
\mathrm{``det}\text{''}
\left(
\begin{array}{cccc}
T_{11} & T_{12} & \dots & T_{1n} \\
T_{21} & T_{22} & \dots & T_{2n} \\
\vdots & & & \vdots \\
T_{n1} & T_{n2} & \dots & T_{nn}
\end{array}
\right) =
\sum_{\substack{\sigma \in P(n) \\ \sigma = \lbrace i_1, \dots , i_n 
\rbrace }} (-1)^{\sgn(\sigma)} T_{1 i_1} T_{2 i_2} \dots T_{n i_n}
\end{equation*}
where the products are given in the order indicated and $ P(n) $ denotes 
the permutations of $ \lbrace 1, \dots , n \rbrace $ and $ \sgn(\sigma) $ 
denotes the sign of the permutation of $ \sigma $.
\end{mydef}

\begin{lem} \label{opdet2}
\begin{multline*}
\mathrm{``det}\text{''} \left (
\begin{array}{cccc}
H_1 & H_{i_1} & \dots & H_{i_p} \\
\underline{f}_1 & \underline{f}_{i_1} & \dots & \underline{f}_{i_p} \\
Q_1^{(1)} & Q_{i_1}^{(1)} & \dots & Q_{i_p}^{(1)} \\
\vdots & & & \vdots \\
Q_1^{(p-1)} & Q_{i_1}^{(p-1)} & \dots & Q_{i_p}^{(p-1)}
\end{array} \right) \\
= p![H_1 \underline{f}_{i_1} Q_{i_2}^{(1)} \dots Q_{i_p}^{(p-1)} + 
\sum_{l=1}^{p} (-1)^l \underline{f}_1 H_{i_l} Q_{i_1}^{(1)} \dots 
Q_{i_{l-1}}^{(l-1)} Q_{i_{l+1}}^{(l)} \dots Q_{i_p}^{(p-1)}] \text{.}
\end{multline*}
\end{lem}

\begin{lem} \label{opdet3}
Suppose $ F(\cdot) $ is $ m \times \infty $ with $ F(z)F(z)^{\ast} $ having 
maximum rank $ p < m $. Suppose that for some function $ \underline{u} : E 
\rightarrow l^2 $, $ F(z)\underline{u}(z) = H(z) $, where $ H = (h_1, \dots 
, h_m) $ with $ h_j \in M(H(E)) $. Then for any $ \pi \in \Pi_{(p+1)}(n) $ 
with $ \pi = (j_1, \dots , j_{p+1}) $, we have for each $ z \in E $
\begin{equation*}
\mathrm{``det}\text{''} \left (
\begin{array}{ccc}
H_{i_1} & \dots & H_{i_{p+1}} \\
\underline{f}_{i_1} & \dots & \underline{f}_{i_{p+1}} \\
Q_{i_1}^{(1)} & \dots & Q_{i_{p+1}}^{(1)} \\
\vdots & & \vdots \\
Q_{i_1}^{(p-1)} & \dots & Q_{i_{p+1}}^{(p-1)}
\end{array} \right) = 0 \text{.}
\end{equation*}
\end{lem}

\section{Proof of Wolff's Theorem for Matrices}

We are now ready to prove our theorem.

\begin{proof}

The solution $ G $ that we seek can be written as a sum of vectors $ G_1, \dots , G_m $. We will find the vector $ G_1 $ here; vectors $ G_2, \dots , G_m $ are found similarly.

Let $ k \leq m $. By the same process we used to obtain \eqref{det}, 
we see that
\begin{equation*}
[\sum_{\substack{\pi \in \Pi_k(m) \\ \pi = \lbrace i_1, \dots i_k \rbrace}} 
\underline{f}_{i_1}(z) Q_{i_2}^{(1)}(z) \dots Q_{i_k}^{(k-1)}(z) Q_{i_k}^{(k-1)\ast}(z) 
\dots Q_{i_2}^{(1)\ast}(z)\underline{f}_{i_1}^{\ast}(z)]^{\frac{3}{2}} \geq \mid 
h_i(z) \mid
\end{equation*}
for all $ z \in \D $ and $ i = 1, \dots , m $. Using 
\eqref{Hinfinitynormestimate}, 
\begin{equation*}
\parallel M_{\underline{f}_{i_1}} M_{Q_{i_2}^{(1)}} \dots M_{Q_{i_k}^{(k-1)}} \parallel 
< \infty \text{,}
\end{equation*}
so by Wolff's Theorem there exists, for each $ \pi \in \Pi_k(m) $, a vector 
$ v_{\pi} $ with entries in $ H^{\infty}(\D) $ such that
\begin{equation*}
k! \sum_{\substack{\pi \in \Pi_k(m) \\ \pi = \lbrace i_1, \dots i_k 
\rbrace}} \underline{f}_{i_1} Q_{i_2}^{(1)} \dots Q_{i_k}^{(k-1)} v_{\pi} = h_1 
\end{equation*}
and
\begin{equation*}
\| M_{v_{\pi}} \| < \infty \text{.}
\end{equation*}

We can rewrite this equation in terms of exterior algebras as
\begin{equation*}
\left(\begin{array}{c} \underline{f}_1 \\ \vdots \\ \underline{f}_m \end{array}\right) 
\wedge
\left(\begin{array}{c} Q_1^{(1)} \\ \vdots \\ Q_m^{(1)} \end{array}\right) 
\wedge \dots \wedge
\left(\begin{array}{c} Q_1^{(k-1)} \\ \vdots \\ Q_m^{(k-1)} 
\end{array}\right) \cdot v_1 = h_1
\end{equation*}
where $ v_1 $ is a vector with $ \left(\begin{array}{c} m \\ k \end{array} 
\right) $ 
entries $ v_{\pi} $ for $ \pi \in \Pi_k(m) $. Then we have $ \| M_{v_i} \| < \infty $.

We claim the vector
\begin{equation*}
G_1 = k\left(\begin{array}{c}1\\0\\ \vdots \\0\end{array}\right) \wedge 
\left(\begin{array}{c} Q_1^{(1)} \\ Q_2^{(1)} \\ \vdots \\ Q_m^{(1)} 
\end{array}\right) \wedge \dots \wedge
\left(\begin{array}{c} Q_1^{(k-1)} \\ Q_2^{(k-1)} \\ \vdots \\ Q_m^{(k-1)} 
\end{array}\right) \cdot v_1
\end{equation*}
is the vector we seek. To prove this, we will consider a more general 
vector,
\begin{equation*}
A = k\left(\begin{array}{c} \alpha_1 I \\ \alpha_2 I \\ \vdots \\ \alpha_m 
I \end{array}\right) \wedge \left(\begin{array}{c} Q_1^{(1)} \\ Q_2^{(1)} 
\\ \vdots \\ Q_m^{(1)} \end{array}\right) \wedge \dots \wedge
\left(\begin{array}{c} Q_1^{(k-1)} \\ Q_2^{(k-1)} \\ \vdots \\ Q_m^{(k-1)} 
\end{array}\right)
\end{equation*}
where $ \alpha_1, \dots , \alpha_m \in H^{\infty}(\D) $. Now
\begin{align*}
\underline{f}_1 A &= k\left(\begin{array}{c} \alpha_1 \underline{f}_1 \\ \vdots \\ \alpha_m \underline{f}_1 
\end{array}\right) \wedge \dots \wedge
\left(\begin{array}{c} Q_1^{(k-1)} \\ \vdots \\ Q_m^{(k-1)} 
\end{array}\right) \\
&= k! \left(\begin{array}{c} \alpha_1 \underline{f}_1 \\ \vdots \\ \alpha_m \underline{f}_1 
\end{array}\right) \wedge 
\sum_{\substack{\sigma \in \Pi_{k-1}(m) \\ \sigma = (j_2, \dots , j_k)}} 
Q_{j_2}^{(1)} \dots Q_{j_k}^{(k-1)} e_{\sigma} \\
&= k! \sum_{j=1}^{m} \sum_{\substack{\sigma \in \Pi_{k-1}(m) \\ \sigma = 
(j_2, \dots , j_k) \\ 1 \notin \sigma}} \alpha_j \underline{f}_1 Q_{j_2}^{(1)} \dots 
Q_{j_k}^{(k-1)} e_j \wedge e_{\sigma} \\
&= k! \alpha_1 \sum_{\substack{\sigma \in \Pi_{k-1}(m) \\ \sigma = (j_2, 
\dots , j_k) \\ 1 \notin \sigma}} \underline{f}_1 Q_{j_2}^{(1)} \dots Q_{j_k}^{(k-1)} 
e_1 \wedge e_{\sigma} \\
&+ k! \sum_{j=2}^{m} \sum_{\substack{\sigma \in \Pi_{k-1}(m) \\ \sigma = 
(j_2, \dots , j_k) \\ 1, j \notin \sigma}} \alpha_j \underline{f}_1 Q_{j_2}^{(1)} \dots 
Q_{j_k}^{(k-1)} e_j \wedge e_{\sigma} \\
&= k! \alpha_1 \sum_{\substack{\pi \in \Pi_{k}(m) \\ \pi = (1, i_2, \dots , 
i_k)}} \underline{f}_1 Q_{i_2}^{(1)} \dots Q_{i_k}^{(k-1)} e_{\pi} \\
&+ k! \sum_{\substack{\pi \in \Pi_{k}(m) \\ \pi = (i_1,\dots , i_k) \\ 1 
\notin \pi}} \sum_{l=1}^{k} (-1)^{l-1} \alpha_{i_l} \underline{f}_1 Q_{i_1}^{(1)} \dots 
Q_{i_{l-1}}^{(l-1)} Q_{i_{l+1}}^{(l)} \dots Q_{i_k}^{(k-1)} e_{\pi} 
\text{.}
\end{align*}
For the second and third equalities, we simply applied the definition of 
the exterior product. For the fourth, we broke the summation into two parts. The 
last inequality is obtained by renaming the indices and using 
\eqref{anticommute}.

By Lemma \ref{opdet2},
\begin{align*}
&k! \sum_{l=1}^{k} (-1)^{l} \alpha_{i_l} \underline{f}_1 Q_{i_1}^{(1)} \dots 
Q_{i_{l-1}}^{(l-1)} Q_{i_{l+1}}^{(l)} \dots Q_{i_k}^{(k-1)} \\
&= \text{``det''} \left(\begin{array}{cccc}
\alpha_1 & \alpha_{i_1} & \dots & \alpha_{i_k} \\ \underline{f}_1 & \underline{f}_{i_1} & \dots & 
\underline{f}_{i_k} \\
\vdots & \vdots & \ddots & \vdots \\
Q_1^{(k-1)} & Q_{i_1}^{(k-1)} & \dots & Q_{i_k}^{(k-1)}
\end{array} \right) 
- k! \alpha_1 \underline{f}_{i_1} Q_{i_2}^{(1)} \dots Q_{i_k}^{(k-1)} \text{.}
\end{align*}
But $ \rank F(z)F(z)^{\ast} = k $, so by Lemma \ref{opdet3}, the $ 
\text{``det''}$ term equals 0. Thus
\begin{align*}
k! \sum_{l=1}^{k} (-1)^{l-1} \alpha_{i_l} \underline{f}_1 Q_{i_1}^{(1)} \dots 
Q_{i_{l-1}}^{(l-1)} Q_{i_{l+1}}^{(l)} \dots Q_{i_k}^{(k-1)} = k! \alpha_1 
\underline{f}_{i_1} Q_{i_2}^{(1)} \dots Q_{i_k}^{(k-1)} \text{.}
\end{align*}
Now we have 
\begin{align*}
\underline{f}_1 A &= k! \alpha_1 \sum_{\substack{\pi \in \Pi_{k}(m) \\ \pi = (1, i_2, 
\dots , i_k)}} \underline{f}_1 Q_{i_2}^{(1)} \dots Q_{i_k}^{(k-1)} e_{\pi} \\
&+ k! \alpha_1 \sum_{\substack{\pi \in \Pi_{k}(m) \\ \pi = (i_1, i_2, \dots 
, i_k) \\ 1 \notin \pi}} \underline{f}_{i_1} Q_{i_2}^{(1)} \dots Q_{i_k}^{(k-1)} 
e_{\pi} \\
&= k! \alpha_1 \sum_{\substack{\pi \in \Pi_{k}(m) \\ \pi = (i_1, i_2, \dots 
, i_k)}} \underline{f}_{i_1} Q_{i_2}^{(1)} \dots Q_{i_k}^{(k-1)} e_{\pi} \\
&= \alpha_1 \left(\begin{array}{c} \underline{f}_1 \\ \vdots \\ \underline{f}_m \end{array}\right) 
\wedge \left(\begin{array}{c} Q_1^{(1)} \\ \vdots \\ Q_m^{(1)} 
\end{array}\right) \wedge \dots \wedge
\left(\begin{array}{c} Q_1^{(k-1)} \\ \vdots \\ Q_m^{(k-1)} 
\end{array}\right) \text{.}
\end{align*}
Thus
\begin{align*}
\underline{f}_1 G_1 &= 1 \cdot \left(\begin{array}{c} \underline{f}_1 \\ \vdots \\ \underline{f}_m 
\end{array}\right) \wedge \left(\begin{array}{c} Q_1^{(1)} \\ \vdots \\ 
Q_m^{(1)} \end{array}\right) \wedge \dots \wedge
\left(\begin{array}{c} Q_1^{(k-1)} \\ \vdots \\ Q_m^{(k-1)} 
\end{array}\right) \cdot v_1 \\
&= h_1
\end{align*}
and similarly,
\begin{align*}
\underline{f}_i G_1 &= 0 \cdot \left(\begin{array}{c} \underline{f}_1 \\ \vdots \\ \underline{f}_m 
\end{array}\right) \wedge \left(\begin{array}{c} Q_1^{(1)} \\ \vdots \\ 
Q_m^{(1)} \end{array}\right) \wedge \dots \wedge
\left(\begin{array}{c} Q_1^{(k-1)} \\ \vdots \\ Q_m^{(k-1)} 
\end{array}\right) \cdot v_1 \\
&= 0
\end{align*}
for $ i \neq 1 $.

We need only show that $ \| M_G \| < \infty $. By 
\eqref{Hinfinitynormestimate},
\begin{align*}
\| M_{G_1} \| & = \| k\left(\begin{array}{c}1\\0\\ \vdots 
\\0\end{array}\right) \wedge 
\left(\begin{array}{c} Q_1^{(1)} \\ Q_2^{(1)} \\ \vdots \\ Q_m^{(1)} 
\end{array}\right) \wedge \dots \wedge
\left(\begin{array}{c} Q_1^{(k-1)} \\ Q_2^{(k-1)} \\ \vdots \\ Q_m^{(k-1)} 
\end{array}\right) \cdot v_1 \| \\
& \leq k \| \left(\begin{array}{c}1\\0\\ \vdots 
\\0\end{array}\right) \wedge 
\left(\begin{array}{c} Q_1^{(1)} \\ Q_2^{(1)} \\ \vdots \\ Q_m^{(1)} 
\end{array}\right) \wedge \dots \wedge
\left(\begin{array}{c} Q_1^{(k-1)} \\ Q_2^{(k-1)} \\ \vdots \\ Q_m^{(k-1)} 
\end{array}\right) \| \| M_{v_1} \| \\
&= k! \| \left(\begin{array}{c}1\\0\\ \vdots \\0 \end{array} \right) \wedge
\sum_{\substack{\pi \in \Pi_{k}(m) \\ \pi = {1, i_2, \dots , i_k}}} 
Q_{i_2}^{(1)} \dots Q_{i_k}^{(k-1)} \|  \| M_{v_1} \| \\
&= k! \| \sum_{\substack{\pi \in \Pi_k(m) \\ \pi = {1, i_2, \dots , i_k}}}
Q_{i_2}^{(1)} \dots Q_{i_k}^{(k-1)} \| \| M_{v_1} \| \\
& \leq k! \left( \sum_{\substack{\pi \in \Pi_k(m) \\ \pi = {1, i_2, \dots , i_k}}} \| M_{Q_{i_2}^{(1)}} \| \dots \| M_{Q_{i_k}^{(k-1)}} \| \right) \|  M_{v_1} \| 
\\ & \leq k! \left( \sum_{\substack{\pi \in \Pi_k(m) \\ \pi = {1, i_2, \dots , i_k}}} \| M_{f_{i_2}} \| \dots \| M_{f_{i_k}} \| \right) \| M_{v_1} \| \\
& \leq k!  \sum_{\substack{\pi \in \Pi_k(m) \\ \pi = {1, i_2, \dots , i_k}}} \| M_{v_1} \| \\ 
& = k! \left( \begin{array}{c} m-1 \\ k-1 
\end{array} \right) \| M_{v_1} \| \text{.}
\end{align*}

Using the estimates in \cite{estimate} for $ \alpha(t) = t^{\frac{1}{2}} $, 
we obtain
\begin{equation*}
\| M_{v_1} \| \leq 1 + 4 \sqrt{e} + 8 \sqrt{2} e + 72 
e^{\frac{3}{2}} = K < 362
\end{equation*}
so
\begin{equation*}
\| M_G \| \leq mk! \left( \begin{array}{c} m-1 \\ k-1 \end{array} \right)
\| M_{v_1} \| \leq mk! \left( \begin{array}{c} m-1 \\ k-1 \end{array} \right) \frac{K}{k!} \leq m \left( \begin{array}{c} m-1 \\ k-1 \end{array} \right) K \text{.}
\end{equation*}
This concludes our proof.

\end{proof}

\section{Further Results}

\subsection{Improved Estimates}

As noted in the introduction, one can improve the estimate in 
Wolff's theorem. The exponent ``$ \frac{3}{2} $'' used in our hypotheses 
isn't optimal, but was used for convenience. 

Cegrell \cite{cegrell} showed that \eqref{eq:scalar} can be replaced with
\begin{equation} \label{alphaestimate}
F(z)F(z)^{\ast} \alpha(F(z)F(z)^{\ast}) \geq |h(z)| \all z \in \D
\end{equation}
where
\begin{equation} \label{alpha}
\alpha(t) = A_0 (\ln \frac{c}{t})^{-\frac{3}{2}}(\ln \ln \frac{c}{t})^{-\frac{3}{2}}(\ln \ln \ln \frac{c}{t})^{-1}
\end{equation}
for $ t \in (0, 1] $ and $ \alpha(0) = 0 $, and for $ F(z) = (f_1(z), \dots , f_n(z)) $. Here $ c $ is chosen so that all log expressions are positive, and $ A_0 $ is chosen so that $ \alpha(1) = 1$.
 Trent \cite{estimate} improved on this estimate (and also allowed for infinitely many functions $ f_i $). The best estimate is currently due to Treil \cite{treil3}. These and further improvements in the estimate automatically carry over to our theorem.

\subsection{Extensions to Other Spaces}

Although we restricted our attention to functions in $ H^{\infty}(\D) $, 
the methods used in the proof of Wolff's theorem for matrices apply to any 
algebra of functions that satisfies a Wolff theorem. (Note that some 
hypotheses may have to be changed. For example, on $ H^{\infty}(\D) $, $ \| 
M_F \| = \| M_F^T \| $, but on other spaces we may have to stipulate that $ 
\max \lbrace \| M_F \| , \| M_F^T\| \rbrace < \infty $.)

As an example, consider Dirichlet space, $ \Dir ^2 (\D) $, defined by
\begin{align*}
\Dir^2(\D) = &\lbrace f : \D \rightarrow \C \mid f \text{ is analytic on } 
\D,&\\ &f(z) = \sum_{n=0}^{\infty} a_n z^n, \| f \| ^2 = 
\sum_{n=0}^{\infty} 
(n+1) | a_n | ^2 < \infty \rbrace \text{.}
\end{align*}
Banjade \cite{banjade} has recently proved that the algebra of multipliers 
on Dirichlet space, $ M(\Dir ^2(\D)) $, satisfies a Wolff theorem: given $ 
\lbrace f_j
\rbrace_{j=1}^{\infty} \subset M(\Dir ^2(\D)) $ and $ h \in M(\Dir ^2(\D)) 
$ such that \eqref{alphaestimate} holds,
\begin{equation*}
|F'(z)F^{\ast}(z)| \alpha(F(z)F(z)^{\ast}) \geq | h'(z) | \all z 
\in \D 
\end{equation*}
and
\begin{equation*}
\| M_F \| \leq 1 \text{,}
\end{equation*}
where $ F(z) = (f_1(z), f_2(z), \dots ) $ and $ \alpha $ is as in 
\eqref{alpha},
then there exists $ \lbrace g_j \rbrace_{j=1}^{\infty} \subset M(\Dir 
^2(\D)) $ with 
$ G(z) = (g_1(z), g_2(z), \dots ) $ such that $ F(z)G(z)^T = h(z) \all z 
\in \D $, and $
\| M_G \| < \infty $.

Thus, given an $ m \times \infty $ matrix $ \F(z) $ of functions in $
 M(\Dir ^2(\D)) $ with $ \max\{\rank  \F(z) \mid z \in \D \} = k \leq m $ 
and an
$ m \times 1 $ vector of functions in $ M(\Dir ^2(\D)) $ such that
\begin{itemize}
\item[(i)] $ \det_k (\F(z)\F(z)^{\ast}) \alpha(\det_k(\F(z)\F(z)^{\ast})) 
\geq |h_i(z)| 
\all z \in \D $, $ i = 1, \dots , m $
\item[(ii)] $ \det_k (| \F'(z)\F(z)^{\ast} | ) \alpha(\det_k 
(\F(z)\F(z)^{\ast})) \geq
| {h'}_i(z) | \all z \in \D $, $ i = 1, \dots , m $
\item[(iii)] $ \| M_{\F} \| = 1 $
\item[(iv)] there exists a function $ \underline{u} : \D \rightarrow l^2 $ 
such that $ \F \underline{u} = H $ everywhere on $ \D $
\end{itemize}
then there exists an $ \infty \times 1 $ vector $ \G(z) $ of functions in $ 
M(\Dir ^2(\D)) $ such that
\begin{itemize}
\item[(a)] $ \F(z)\G(z)^T = H(z) \all z \in \D, $ and
\item[(b)] $ \| M_{\G} \| < \infty $.
\end{itemize}
Note that in this and other cases outside of $ H^{\infty}(\D) $ we would use Lemma 1 from \cite{mctII} instead of \eqref{Hinfinitynormestimate} to 
estimate $ \| M_G \| $.

\subsection{Radicals}

As previously noted, Wolff's condition \eqref{eq:wolff} is not sufficient 
to show $h \in \I $. However, it 
is necessary \emph{and} sufficient to show that $ h $ is contained in the 
\emph{radical} of $ \I $.

We would like to show a similar result for the matrix case. Let $ F $ and $ 
H $ be as before, with $ \det_k (F(z)F(z)^{\ast}) \geq | h_i(z) |^n \all z 
\in \D, \; i = 1, \dots , m $, and for some $ n \in \mathbb{N} $, where $ k = 
\max \lbrace \rank F(z) | z \in \D \rbrace $. Suppose also 
that $ \| M_F \| = 1 $ and that we can find a $ \underline{u} $ such that $ 
F\underline{u} = H $ on $ \D $, as before. Then by Wolff's Theorem for Matrices, there 
exists an $ \infty \times 1 $-vector $ G $ with entries in $ H^{\infty}(\D) 
$ such that $ FG = H^{3n} $ everywhere on $ \D $. (Here $ H^n $ is the vector 
obtained by raising each entry of $ H $ to the $ n $th power.) 

On the other hand, suppose we have $ FG = H^n $ for some $ n \in \mathbb{N} 
$. Then
\begin{align*}
FG(G^{\ast}F^{\ast}) &= H^n (H^n)^{\ast} \Rightarrow \\
F \| M_G \|^2 F^{\ast} & \geq H^n (H^n)^{\ast} \Rightarrow \\
{\det}_1( \| M_G \|^2 FF^{\ast} ) & \geq {\det}_1(H^n(H^n)^{\ast})
= \sum_{i=1}^{m} |h_i|^{2n} \Rightarrow \\
C \cdot{\det}_1(FF^{\ast}) & \geq |h_i|^{2n}
\end{align*}
$ \all z \in \D, \; i = 1, \dots , m $, where $ C =  \|M_G\|^{2m} $.

Note that if $ k = \max \lbrace \rank F(z) | n \in \D \rbrace = 1 $, then 
this second statement is the converse of the first. It is currently unknown 
whether the converse holds for $ k > 1 $. 

\subsection{When $H$ Is a Matrix}

Our theorem extends easily to the case where $ H $ is an $ m \times n $ 
matrix. If $ F $ is $ m \times \infty $, we seek an $ \infty \times n $ 
matrix $ G $ such that $ FG = H $. Wolff's Theorem for Matrices allows us 
to find $ G $ by finding its $ n $ columns $ g_1, \dots , g_n $.

What if we wish to solve an equation involving two (or more) matrices $ F_1 
$ and $ F_2 $? That is, we wish to find $ G_1 $ and $ G_2 $ such that
\begin{equation*}
F_1 G_1 + F_2 G_2 = H \text{.}
\end{equation*}
This is handled easily if we define $ \mathcal{F} = [ F_1 \, F_2 ] $; that 
is, $ \mathcal{F} $ is obtained by concatenating $ F_1 $ with $ F_2 $ 
(rearranging the entries in the case where $ F_1 $ and $ F_2 $ are $ m 
\times \infty $). Then, provided hypotheses \emph{(i)}, \emph{(ii)}, and
\emph{(iii)} of our main theorem hold on $ \mathcal{F} $,
there exists $ \mathcal{G} $ such that
\begin{equation*}
\mathcal{F} \mathcal{G} = H \text{,}
\end{equation*}
and from $ \mathcal{G} $ we obtain $ G_1 $ and $ G_2 $ such that
\begin{equation*}
\mathcal{F} \mathcal{G} = F_1 G_1 + F_2 G_2 \text{.}
\end{equation*}

\end{document}